\documentclass{article}
\usepackage{amssymb}
\newtheorem{thm}{Theorem} 
\newtheorem{cor}[thm]{Corollary}

\newtheorem{prop}[thm]{Proposition}
\newtheorem{rem}[thm]{Remark}
\newcommand{\wt}{\beta}
\newcommand{\One}{\vec{\mathbf{1}}}

\begin{document}
\title{A Note on Diffusion State Distance}  
\author{Neal Madras \\ Department of Mathematics and Statistics \\
York University \\ 4700 Keele Street  \\ Toronto, Ontario  M3J 1P3 Canada 
\\  {\tt  madras@mathstat.yorku.ca} }
\maketitle

\begin{abstract}
Diffusion state distance (DSD) is a metric on the vertices of
a graph, motivated by bioinformatic modeling.  
Previous results on the convergence of DSD to a limiting metric 
relied on the definition being based on symmetric or reversible
random walk on the graph.  We show that convergence holds
even when the DSD is based on general finite irreducible Markov chains.
The proofs rely on classical potential theory of Kemeny and Snell.
\end{abstract}

\section{Introduction and Main Results}
   \label{sec-intro}
In 2013, Cao \textit{et al.}\ \cite{Cao} defined a metric (or rather,
a sequence of metrics) on the set of vertices of a finite undirected 
graph, motivated by functionality considerations in protein interaction 
networks.  
In this class of metrics, called ``diffusion state distance'' (DSD),
\cite{Cao} considered a random walk on the vertices of the graph,
and assessed
the closeness of two states $u$ and $v$ by comparing 
the expected number of visits to all states (within a given time 
horizon) when the initial state 
is $u$ and when the initial state is $v$.
(See below for the mathematical details.)
The following year,  Cao \textit{et al.}\ \cite{Cao2} extended the
metric to graphs with weighted edges (representing confidence in
the presence of the edge), which leads to consideration of
Markov chains instead of random walks.  
The paper  \cite{Cao2} also permitted directed edges.
Further generalizations were considered in \cite{Boe}, in 
particular to bipartite graphs (which lead to periodic random 
walks).

To proceed, we now describe  a general framework for random walks 
and Markov chains on graphs. 
Let $V$ be  a (finite) set of points, representing the
vertices or nodes of our network. 
For each (ordered) pair of vertices $i$ and $j$, let $\wt_{ij}$ be a 
nonnegative number associated with the directed edge from $i$ to $j$.
(In \cite{Cao2}, the magnitude of $\wt_{ij}$ corresponds to the degree
of certainty that this edge is present. If there is certainly no 
edge from $i$ to $j$, then $\wt_{ij}$ is zero.)
We consider the Markov chain on $V$ with the property that 
when the chain is at a node $i$, it chooses its next node with 
probability proportional to the weights on the edges leading out of 
$i$.  That is, the one-step transition probabilities for the 
Markov chain are 
\begin{equation}
   \label{eq.pdef}
      p_{ij} \;=\; \frac{\wt_{ij}}{\sum_{k\in V}\wt_{ik}} 
    \hspace{5mm} \mbox{for $i,j\in V$}. 
\end{equation}
(We assume that there is no $i$ for which 
the denominator in Equation (\ref{eq.pdef}) is 0.)
The case of a classical random walk on a directed graph (with no weights)
is obtained by requiring each $\wt_{ij}$ to equal 0 or 1, according to 
whether or not there is an edge from $i$ to $j$.  The case of
a weighted undirected graph is obtained
by requiring $\wt_{ij}=\wt_{ji}$  for every $i$ and $j$ in $V$; this is
the case assumed by Chung and Yao \cite{CY}, which 
in turn was built upon by \cite{Boe}.

Note that the form of Equation (\ref{eq.pdef}) is completely general,
since the transition probabilities $p_{ij}$ of any given Markov chain 
satisfy (\ref{eq.pdef}) if we let $\wt_{ij}$ be $p_{ij}$ for every $i$ and $j$
(since $\sum_k p_{ik}=1$ in any Markov chain).

We shall denote the set of states of the Markov chain by $V=\{1,\ldots,n\}$.
We shall write $P$ for the (one-step) transition probability matrix
with entries $p_{ij}$.  For $l=0,1,2,\ldots$, we write $p_{ij}^{(l)}$
for the $l$-step transition probabilities, which are the entries
of the matrix $P^l$.
We say that the Markov chain (or $P$) is irreducible if it is possible
to get from every state to every other state; that is, for every 
pair of states $i$ and $j$, there is an $l\geq 0$ such that $p_{ij}^{(l)}>0$.
We say that an irreducible Markov chain (or $P$) is periodic if there 
is an integer $d\geq 2$ such that $p_{11}^{(l)}=0$ whenever $l$ is not
a multiple of $d$; in this case $d$ is called the period.  We say 
that the chain is aperiodic if it is not periodic.

Throughout this paper, except for the generalizations discussed in 
Section \ref{sec-gen},
 \textit{we shall assume that $P$ is the 
transition probability matrix of a finite irreducible Markov chain.} 

\subsection{The Aperiodic Case}
   \label{sec-aper}
In this section, we shall assume that $P$ is
\textit{aperiodic} as well as  finite and irreducible.

Let $\pi$ be the equilibrium distribution of $P$, i.e.\ the 
row vector $(\pi_1,\ldots,\pi_n)$ such that $\pi P =\pi$ and $\sum \pi_i =1$.
(This is written $\pi^T$ in \cite{Cao} and $\alpha$ in \cite{KSK}).  
Under our assumptions, $\pi$ exists and is unique.  Moreover, the limit matrix
\begin{equation}
   \label{eq.W}
     W   \; :=  \;  \lim_{k\rightarrow\infty}P^k   
\end{equation}
exists, and $W_{ij} =\pi_j$ for all states $i$ and $j$.
(This is called $W$ in \cite{Cao}, and $A$ in \cite{KS} and \cite{KSK}.)
Writing $\One$  to denote the column vector with all 
entries equal to 1, we have $W\,=\,\One \pi$.

\smallskip
For states $u$ and $v$ and any integer $k\geq 0$, we follow \cite{Cao} and
define
\[    He^{(k)}(u,v)   \;=\;   \sum_{l=0}^k  p_{uv}^{(l)}   \,.
\]
Thus $He^{(k)}(u,v)$  is the expected number of visits to $v$ 
within the first $k$ steps of the chain,
given that the chain starts at $u$.  (This is called $N^{(k)}_{uv}$ in 
\cite{KSK}, and $\textbf{M}_u[\bar{\textbf{y}}^{(k+1)}_v]$ in \cite{KS}.)

\smallskip
\begin{prop}
   \label{propHe}
For an aperiodic finite irreducible Markov chain:
\\
(a)
The matrix $Z$ defined by
\[    Z \;=\; \sum_{k=0}^{\infty}(P-W)^k  \]
exists and equals $(I-P+W)^{-1}$.
(\textit{Note}:  
\cite{Cao} writes $D$ for $P-W$).
\\
(b)
For all states $u$, $v$, and $w$, we have
\[    \lim_{k\rightarrow\infty} He^{(k)}(u,w)-He^{(k)}(v,w) 
   \hbox{ exists and equals }
    Z_{uw}-Z_{vw}\,.
 \]
\end{prop}

\medskip
\noindent
\begin{rem}
   \label{rem-bZ}
Part (b) agrees with Lemma 3 of \cite{Cao}, since 
$Z_{uw}-Z_{vw}\,=\,(b^T_u-b^T_v)Z\,b_w$ (where $b_i$ is the column 
vector having
$i^{th}$ entry equal to  1 and all other entries equal to 0.)
\end{rem}
 
We now define the diffusion state distance metrics.
For states $u$ and $v$ and integers $k\geq 0$, let 
\begin{equation}
  \label{eq-defk}
   DSD^{(k)}(u,v) \;=\;  \sum_{w=1}^n |He^{(k)}(u,w)-He^{(k)}(v,w)| \,.   
\end{equation}
Lemma 1 of \cite{Cao} shows that $DSD^{(k)}$ is a metric for aperiodic irreducible
random walks; see Section \ref{sec-gen} below for discussion and generalization.
We also write
\begin{equation}
  \label{eq-DSDlim}
     DSD^{(\infty)}(u,v)  \;=\; \lim_{k\rightarrow\infty} DSD^{(k)}(u,v)   
\end{equation}
if this limit exists.

\begin{cor}
  \label{corDSD}
Consider an aperiodic finite irreducible Markov chain.
Then the limit (\ref{eq-DSDlim}) exists, and 
\begin{equation}
   \label{eq-DSD}
      DSD^{(\infty)}(u,v)  \;=\; \sum_w |Z_{uw}-Z_{vw}|  \,.
\end{equation}
Moreover, $DSD^{(\infty)}$  is a metric.
\end{cor}
Corollary \ref{corDSD} was proven in \cite{Cao} for undirected 
unweighted graphs ($\wt_{ij}=\wt_{ji}\in\{0,1\}$), 
and in \cite{Boe} for graphs with symmetric
weights ($\wt_{ij}=\wt_{ji}$).    
(Observe that the Green's function matrix $\mathbb{G}$ in \cite{Boe}
is precisely our $Z$.)
Both proofs rely on the diagonalization of the matrix $P$, which holds
because the symmetry of weights  implies that $P$ can be expressed
as a self-adjoint
operator  (equivalently, that the Markov chain is reversible).
In contrast, the proofs in the present paper rely on the classical
potential theory for general finite (and countable) Markov chains
in \cite{KS} and/or \cite{KSK}.

We remark that \cite{Boe} also considers the $DSD^{(k)}_q$ metric
for $q\geq 1$, 
defined by replacing the $l_1$ norm on the right hand side of 
Equation (\ref{eq-defk}) by the $l_q$ norm, i.e.
\[   DSD^{(k)}_q(u,v) \;=\;
  \left(\sum_{w=1}^n |He^{(k)}(u,w)-He^{(k)}(v,w)|^q\right)^{1/q} \,.   
\]
The extension of Corollary \ref{corDSD} and
Equation (\ref{eq-DSDlim}) 
to $DSD_q$ is immediate.

\subsection{The Periodic Case}
   \label{sec-per}

In this section, we shall consider the case  that $P$ is 
\textit{periodic} as well as finite and irreducible.
An important example  is the case of a random walk 
on a bipartite graph.  This example was considered in \cite{Boe}. 

As explained in Chapter 5.1 of \cite{KS},
the relevant theory for aperiodic Markov chains extends to the periodic case
with some modifications.  There is a unique equilibrium probability
distribution $\pi$ such that $\pi P=\pi$.   We still define the matrix $W$ 
by $W_{ij}=\pi_j$.  However, the limit of Equation (\ref{eq.W}) does not exist
in the usual sense, but the equation is correct for the Cesaro limit.
Moreover, the matrix $I-P+W$ is invertible, and Proposition \ref{propHe}(a)
holds if we interpret the infinite sum defining $Z$ to be the Cesaro sum. 

For now, consider the general case
that $P$ is finite and irreducible (not necessarily periodic).  
For $0\leq \alpha<1$, let $P_{\alpha}=\alpha I+(1-\alpha)P$ be the 
``lazy Markov chain'', which moves according to $P$ except that with probability 
$\alpha$ it decides to stay wherever it is for one time unit. 
Observe that $\pi = \pi P_{\alpha}$, so the 
equilibrium distribution is independent of $\alpha$, as is the matrix $W$.
Let $Z_{\alpha}=(I-P_{\alpha}+W)^{-1}$.
Since this is well-defined, we can use Equation (\ref{eq-DSD})
as our definition of $DSD^{(\infty)}$ for our periodic chain $P$,
even though the limit of Proposition \ref{propHe} may not exist
except in the Cesaro sense.
The following result shows that this definition is consistent
with the usual definition.
Observe that
if $0<\alpha<1$, then $P_{\alpha}$ is aperiodic (and irreducible), 
and so the results of Section \ref{sec-aper}  apply.

\begin{prop}
   \label{prop-period}
Assume that $P$ is finite and irreducible.
\\
(a) For every $\alpha\in [0,1)$, we have 
$Z_{\alpha} \,=\, (1-\alpha)^{-1}(Z_0-\alpha W)$.
\\
(b)  Let $DSD^{(\infty;\alpha)}$
be the $DSD^{(\infty)}$ metric of the Markov chain  $P_{\alpha}$.  Then 
for all states $u$ and $v$,
\[   DSD^{(\infty;\alpha)}(u,v)  \;=\;
     (1-\alpha)^{-1} DSD^{(\infty;0)}(u,v)\,.   \]
\end{prop}
In particular, part (b) generalizes Theorem 1 of \cite{Boe} to 
general Markov chains.  It shows that for periodic chains,
the formal definition of Equation (\ref{eq-DSD}) is consistent
by continuity with the (more directly  motivated) definition (\ref{eq-DSDlim}) for 
aperiodic chains, i.e.\
\[   \lim_{\alpha\rightarrow 0}DSD^{(\infty;\alpha)}(u,v)  \;=\;
     DSD^{(\infty;0)}(u,v)\hspace{5mm}\mbox{for all $u,v$}.   \]

\section{Proofs of Results}
   \label{sec-pfA}

\noindent
\textbf{Proof of Proposition \ref{propHe}}: (a)
This is direct from Propositions 9-75 and 9-76(4) of \cite{KSK},
or Theorem 4.3.1 of \cite{KS}.

\smallskip
\noindent
(b)  This follows from Corollary 4.3.5 of \cite{KS}.  Alternatively,  here is
a proof using the theory for denumerable chains (see discussion at the 
beginning of Section \ref{sec-gen} below).
Definition 9-24 of \cite{KSK} 
defines the matrix $C$ by
\[    C_{ij}   \;=\; \lim_{k\rightarrow\infty}(N_{jj}^{(k)}-N_{ij}^{(k)}) 
  \hspace{3mm} \left(   \;=\;  \lim_{k\rightarrow\infty}
     (He^{(k)}(j,j)-He^{(k)}(i,j) ) \right) \]
whenever this limit exists.  Proposition 9-77\footnote{To explain the 
 notation of \cite{KSK} in this proposition:  $E$ is the matrix of all 1's, 
and $Z_{dg}$ is the diagonal matrix obtained from $Z$ 
 by making all off-diagonal entries 0.} 
of \cite{KSK} shows that $C$ exists for aperiodic finite irreducible 
chains and that
\[    C_{ij} \;=\; Z_{jj}\,-\,Z_{ij}  \,. \]
It follows that 
\begin{eqnarray*}
     \lim_{k\rightarrow\infty} He^{(k)}(u,w)-He^{(k)}(v,w)   
     & = & C_{vw}-C_{uw}   \\
        & = &   (Z_{ww}-Z_{vw}) \,-\,(Z_{ww}-Z_{uw})   \\
        & = & Z_{uw}-Z_{vw}  \,.
\end{eqnarray*}
\hfill $\Box$

\medskip
\noindent
\textbf{Proof of Corollary \ref{corDSD}}:  
The second sentence of the corollary follows immediately from 
Proposition \ref{propHe}.  To show the final sentence, observe that  
Equation (\ref{eq-DSD}) says that $DSD^{(\infty)}(u,v)$ is the $l_1$ distance 
between rows $u$ and $v$ of the matrix $Z$.  Since $Z$ is invertible
(by Proposition \ref{propHe}(a)), it must have distinct rows.  
It follows immediately that $DSD^{(\infty)}$ is a metric.
\hfill $\Box$

\medskip
\noindent
\textbf{Proof of Proposition \ref{prop-period}}:  
By Theorem 5.1.3 of \cite{KS}, we have $Z\One =\One$ and $\pi Z =\pi$,
where $Z$ denotes either $Z_{\alpha}$ or $Z_0$.  
Therefore
\begin{eqnarray}
  \label{eq.ZW}
    ZW  & = & Z\One \pi  \;=\; \One \pi \;=\; W   \mbox{\hspace{5mm}and}  \\
  \label{eq.WZ}
    WZ & = & \One \pi Z  \;=\; \One \pi  \;=\; W \,.
\end{eqnarray} 
Therefore we have
\begin{eqnarray*}
   I \;=\; Z_{\alpha}(I-P_{\alpha}+W)  & = & 
       Z_{\alpha}\left((1-\alpha)(I-P+W)+\alpha W\right)   \\
    & = & (1-\alpha)Z_{\alpha}Z_0^{-1}\,+\, \alpha Z_{\alpha}W   \\
   & = & (1-\alpha)Z_{\alpha}Z_0^{-1}\,+\, \alpha W
     \hspace{10mm}\mbox{[by Eq.\ (\ref{eq.ZW})]}  \,,
\end{eqnarray*}
which yields   
\[
   Z_0 \;=\; (1-\alpha)Z_{\alpha} \,+\,\alpha WZ_0  
   =\; (1-\alpha)Z_{\alpha} \,+\,\alpha W 
        \hspace{10mm}\mbox{[by Eq.\ (\ref{eq.WZ})]}.\]
Part (a) follows immediately. 

Next, recall the observation and notation of 
Remark \ref{rem-bZ} above.  By this and Equation (\ref{eq-DSD}), 
we see that part (b) will follow if we can prove that 
$\gamma Z_{\alpha} \,=\, (1-\alpha)^{-1}\gamma Z_0$, where
$\gamma=b^T_u-b^T_v$.
But this is a direct consequence of part (a) and the fact that
$\gamma W \,=\, \gamma \One \pi \,=\, 0$.
\hfill $\Box$

\section{Some Generalizations}
   \label{sec-gen}

Since much of the potential theory for finite Markov chains in \cite{KS}
also extends to countable chains \cite{KSK}, it is of theoretical 
interest to consider whether DSD theory extends to chains with 
countably many states.  
Our proof of Proposition \ref{propHe} using results of \cite{KSK} shows
that this proposition extends to 
Markov chains that are \textit{strong ergodic} (i.e.\ chains that are
positive recurrent and have the property that
for all states $i$ and $j$, 
the expected square of the time to reach $j$ [when starting from $i$] is 
finite);  in  \cite{KSK}, see Section  9.5,
p.\ 274, as well as Definition 9-71.
However it is less clear whether Corollary \ref{corDSD} holds
in this case, since even the convergence of the infinite sum is
not obvious.

Finally, we show that $DSD^{(k)}$ is a metric 
for every finite or countable Markov chain except for some very special
cases.  For one counterexample, consider
the irreducible Markov chain with states $V=\{1,2,3,4\}$ such that
$p_{12}=p_{13}=0.5$, 
$p_{24}=p_{34}=p_{41}=1$,   
and $p_{ij}=0$ otherwise.  It is not hard to check that
$He^{(2)}(1,i)=He^{(2)}(4,i)$ for every $i$, and so 
$DSD^{(2)}$ is not a metric.  
Note that this chain has period $3$.
The next result shows that all counterexamples have features in 
common with this example.

\begin{prop}
  \label{prop-metric}
Let $P$ be the matrix of a (finite or countable) Markov chain, 
and let $k$ be a finite nonnegative integer.  Assume that 
there do not exist two distinct states $i$ and $j$ such that 
$p^{(k+1)}_{ii}=p^{(k+1)}_{jj}=1$ and $p^{(t)}_{ij}>0$ for some $t$.
Then $DSD^{(k)}$ is a metric.
\end{prop}
In particular, if $P$ is aperiodic or if the period of $P$ does not
divide $k+1$, then $DSD^{(k)}$ is a metric.  But the condition
$p^{(k+1)}_{ii}=1$ is much more restrictive than requiring
the period to divide $k+1$.
Also note that Proposition \ref{prop-metric} does not assume
that $P$ is irreducible.

The following proof is essentially a streamlined generalization of the 
proof of Lemma 1 in \cite{Cao}.

\smallskip
\noindent
\textbf{Proof of Proposition \ref{prop-metric}}:
Write $N^{(t)}=1+P+\cdots+P^t$ for $t=0,1,\ldots$.
Since $DSD^{(k)}(u,v)$ is the the $l_1$ distance between rows
$u$ and $v$ of the (possibly infinite) matrix $N^{(k)}$, it 
suffices to show that all rows of $N^{(k)}$ are distinct.

Assume that $u$ and $v$ are distinct states such that
rows $u$ and $v$ of $N^{(k)}$ are the same.  For vectors $b_i$
as defined in Remark \ref{rem-bZ}, let $\gamma = b^T_u-b^T_v$.
Then $\gamma N^{(k)}=0$.  Hence $\gamma N^{(k)}P=0$.  From the definition,
we have $N^{(k)}P +I\,=\,N^{(k)}+P^{k+1}$, and hence it follows
that 
\[   \gamma = \gamma P^{k+1}.   \]  
Multiplying this equation on the right by $b_u$ gives the equation
$1=(P^{k+1})_{uu}-(P^{k+1})_{vu}$.  This can only happen if 
$(P^{k+1})_{uu}=1$.  Similarly, we deduce $(P^{k+1})_{vv}=1$.
Finally, since $N^{(k)}_{uv}=N^{(k)}_{vv}$ and $N^{(k)}_{vv}\geq I_{vv}=1$,
there must be a $t\leq k$ such that $(P^{t})_{uv}>0$.
This contradicts the hypothesis of the proposition, and thus the proof
is complete.
\hfill $\Box$

\section{Acknowledgments}
I am grateful to Lenore Cowen for helpful conversations and her enthusiasm.
This research was supported in part by a Discovery Grant from NSERC Canada.
Part of this work was done while visiting the Fields Institute for Research in 
Mathematical Sciences.


\begin{thebibliography}{999}

\bibitem{Boe}  Boehnlein, E., Chin, P., Sinha, A., and Lu, L. (2014), ``Computing Diffusion
State Distance using Green's Function and Heat Kernel on Graphs.'' In  \textit{Algorithms and Models
for the Web Graph}, Bonato, A., Graham, F.C., and Pralat, P. (eds.), 
\textit{Lecture Notes in Computer Science}  \textbf{8882}, 79--95.  (Also at
arXiv:1410.3168v1.)


\bibitem{Cao}  Cao, M. \textit{et al.} (2013), ``Going the Distance 
for Protein Function Prediction:  A New Distance Metric for Protein 
Interaction Networks'', \textit{PLoS ONE} \textbf{8}(10): e76339.

\bibitem{Cao2}  Cao, M. \textit{et al.} (2014), ``New Directions for
Diffusion-Based Network Prediction of Protein Function:  Incorporating
Pathways with Confidence,'' \textit{Bioinformatics} \textbf{30}, i219--i227;
doi:10.1093/bioinformatics/btu263. 

\bibitem{CY} Chung, F. and Yao, S.-T. (2000), ``Discrete Green's functions,''
J. Combinatorial Theory A, \textbf{91}, 191--214.

\bibitem{KS}  Kemeny, J.G., and  Snell, J.L. (1976), 
\textit{Finite Markov Chains}. Springer, New York.

\bibitem{KSK}  Kemeny, J.G.,  Snell, J.L.,  and Knapp, A.W. (1976), 
\textit{Denumerable Markov Chains}, Springer.

\end{thebibliography}
\end{document}